\documentclass[11pt]{article}
\input epsf
\usepackage{amsfonts}
\usepackage{amscd}
\usepackage{amsmath,amssymb}
\usepackage{amstext}
\usepackage{amscd}
\usepackage[all]{xy}
\usepackage{pstricks,pst-node,pst-tree}

\newtheorem{lem}{LEMMA}[section]
\newtheorem{theo}[lem]{THEOREM}

\newtheorem{coro}[lem]{COROLLARY}
\newtheorem{prop}[lem]{PROPOSITION}

\newtheorem{definition}[lem]{Definition}

\newtheorem{rem}[lem]{Remark}

\renewcommand{\descriptionlabel}[1]%
       {\hspace{\labelsep}\textsf{#1}}

\begin{document}

\title{Any smooth knot $\mathbb{S}^{n}\hookrightarrow\mathbb{R}^{n+2}$ is isotopic to a cubic knot contained in the canonical scaffolding of $\mathbb{R}^{n+2}$
}
\author{Margareta Boege\thanks{Research partially supported by PFAMU-DGAPA.}, Gabriela Hinojosa\thanks{Research  partially supported by CONACyT CB-2007/83885
and Promep-Pifi3.3 CA.},
Alberto Verjovsky\thanks{Research partially
supported by CONACyT project U1 55084, and PAPIIT (Universidad
Nacional Aut\'onoma de M\'exico) \# IN102108.}}
\date{November 17, 2009}

\maketitle
\rightline{\footnotesize Dedicated to Jos\'e Mar\'{\i}a Montesinos}
\rightline{\footnotesize on the occasion of his $65^{th}$ anniversary}

\begin{abstract}
The $n$-skeleton of the canonical cubulation $\cal C$ of $\mathbb{R}^{n+2}$ into unit cubes is called the {\it canonical scaffolding} ${\cal{S}}$.
In this paper, we prove that any smooth, compact, closed, $n$-dimensional submanifold of $\mathbb{R}^{n+2}$ with trivial normal bundle
can be continuously isotoped by an ambient isotopy to a cubic submanifold  contained in ${\cal{S}}$. In particular,
any smooth knot $\mathbb{S}^{n}\hookrightarrow\mathbb{R}^{n+2}$
can be continuously isotoped to a knot contained in ${\cal{S}}$.

\end{abstract}

\section{Introduction}

In this paper we consider smooth higher dimensional knots, that is, spheres $\Bbb S^n$ smoothly embedded in $\Bbb R^{n+2}$.
In  $\Bbb R^{n+2}$ we have the canonical cubulation  $\cal C$ by translates of the unit $(n+2)$-dimensional cube.
We will call the $n$-skeleton ${\cal{S}}$ of this cubulation the {\it canonical scaffolding} of $\Bbb R^{n+2}$
(see section 2 for precise definitions). We consider the question of whether
it is possible to continuously deform the smooth knot by an ambient isotopy so that the deformed knot is contained in the scaffolding. In particular,
a positive answer to this question implies  that knots can be embedded as cubic sub-complexes of $\Bbb R^{n+2}$, which in turn implies
the well-known fact that smooth knots can be triangulated by a PL triangulation (\cite{cairns}). The problem of embedding an abstract cubic complex into some skeleton
of the canonical cubulation can be traced back to S.P. Novikov. A considerable amount of work has been done regarding this problem (see, for example
\cite{dolbilin}). The question is non-trivial; for instance, among cubic manifolds there are non-combinatorial ones and therefore
non-smoothable ones. Also there is a series of very interesting papers by Louis Funar regarding cubulations of manifolds (\cite{funar}, \cite{funar1}).
The possibility of considering a knot as a cubic submanifold contained in the $n$-skeleton of the canonical cubulation of $\Bbb R^{n+2}$ has many advantages. For instance,
in the important case of classical knots $n=1$,
Matveev and Polyak  \cite{matveev} begin the exposition of finite type invariants from the ``cubic'' point
of view and show how one can clearly describe invariants such as polynomial invariants, Vassiliev-Goussarov
invariants and finite type invariants of three-dimensional integer homology spheres (in this regard see also the
unpublished important paper by Fenn, Rourke and Sanderson \cite{FRS}).
Cubic complexes may play a role in extending these invariants to higher dimensional knots.\\

In this paper we prove that any smooth, compact, closed, $n$-dimensional submanifold of $\Bbb R^{n+2}$ with trivial normal bundle
can be continuously isotoped by an ambient isotopy of $\Bbb R^{n+2}$ onto a cubic submanifold contained in ${\cal{S}}$.
In particular, any knot can be isotoped onto a cubic knot contained in ${\cal{S}}$.

\section{Cubulations for $\mathbb{R}^{n+2}$}

A {\it cubulation} of $\mathbb{R}^{n+2}$ is a decomposition of $\mathbb{R}^{n+2}$ into a collection $\cal C$ of $(n+2)$-dimensional
cubes such that any two of its hypercubes are
either disjoint or meet in one common face of some dimension. This provides  $\mathbb{R}^{n+2}$ with the structure of a cubic
complex.\\

In general, the category of cubic complexes and cubic maps is similar to the simplicial category. The only
difference consists in considering cubes of different dimensions instead of simplexes. In this context,
a cubulation of a manifold is specified by a cubical complex PL homeomorphic to the
manifold (see \cite{dolbilin}, \cite{funar}, \cite{matveev}). \\

The {\it canonical cubulation} $\cal C$ of $\mathbb{R}^{n+2}$ is the decomposition into hypercubes which are the images of the unit cube
$I^{n+2}=\{(x_{1},\ldots,x_{n+2})\,|\,0\leq x_{i}\leq 1\}$ by translations by vectors with integer coefficients.

Consider the homothetic transformation $\frak{h}_{m}:\mathbb{R}^{n+2}\rightarrow\mathbb{R}^{n+2}$ given by
$\frak{h}_m(x)=\frac{1}{m} x$,
where $m>1$ is an integer. The set $\frak{h}_{m}(\cal{C})$ is called a {\it subcubulation} or {\it cubical subdivision} of $\cal C$.

\begin{definition} The $n$-skeleton of $\cal C$, denoted by $\cal S$, consists of the union of the $n$-skeletons of the cubes in  $\cal C$,
{\it i.e.,} the union of all cubes of dimension $n$ contained in the faces of the $(n+2)$-cubes in  $\cal C$.
We will call $\cal S$ the {\it canonical scaffolding} of $\mathbb{R}^{n+2}$.
\end{definition}

Any cubulation of $\mathbb{R}^{n+2}$ is obtained by applying a conformal transformation
$x\mapsto{\lambda{A(x)}+a},\,\,\,\lambda\neq0,\,\,a\in \mathbb{R}^{n+2},\,\,\,A\in{SO(n+2)}$ to the canonical cubulation. \\

In this section, we will prove that the union of the cubes of the canonical cubulation which intersect a fixed hyperplane is a closed tubular
neighborhood of the hyperplane ({\it i.e.,} a bicollar). Hence the boundary of the union of these cubes
has two connected components and therefore the hyperplane can be isotoped to any of them. This isotopy can be realized using normal segments to
the hyperplane. The boundaries of this bicollar are contained in the $(n+1)$-skeleton of the canonical cubulation.
If $P$ is a hyperplane of $\mathbb{R}^{n+2}$ and $\Pi\subset P$ is an $n$-dimensional affine subspace we use the same sort of ideas
to show that the couple $(P,\Pi)$ can be isotopically deformed to the boundary of the bicollar above in such a way that $\Pi$
is deformed into the $n$-skeleton of the canonical cubulation.\\

These linear cases already contain the main ingredients of our proofs which in particular include convexity properties.\\

Later we will generalize these ideas to the case of any smooth, compact, closed
codimension one submanifold $M$ of $\mathbb{R}^{n+2}$: the union of the cubes of a small cubulation ({\it i.e.,} with cubes of sufficiently
small diameter) which intersect $M$ is a bicollar neighborhood of $M$. We can deform $M$ to any of the boundary components using an adapted flow
which is nonsingular in the tubular neighborhood and is transverse to $M$.
Finally using the same type of ideas we prove at the end of the section
that any smooth codimension two submanifold $N$ can be deformed into the $n$-skeleton of a sufficiently small cubulation.

\begin{prop}
Let $P$ be a hyperplane in $\mathbb{R}^{n+2}$. Let $\cal{C}$ be the canonical cubulation of $\mathbb{R}^{n+2}$ and
${\cal {Q}}_{P}=\cup\{Q\in {\cal{C}}:Q\cap P\neq\emptyset\}$. Then:
\begin{enumerate}
\item  $P$ is contained in the interior of ${\cal {Q}}_{P}$.
\item Let $Q\in {\cal{C}}$
such that the distance $d(Q,P)=m>0$. Then the set
${\cal{M}}:=\{q\in{Q}\,\,\,|\,\,\,d(q,P)=m\}$ is a cubic simplex, of some dimension, of the boundary of $Q$.
\item $\partial {\cal {Q}}_{P}$ is the union of $(n+1)$-dimensional cubes which are faces of cubes in ${\cal {Q}}_{P}$.
\end{enumerate}
\end{prop}

{\it Proof.} \\
{\it 1.} Let $p\in {\mathbb R}^{n+2},\,\,n\geq1$ and ${\cal{C}}$ be the canonical cubulation of ${\mathbb R}^{n+2}$. Let us consider the set
$C_{p}=\cup\{Q\in {\cal{C}}:p\in Q\}$.
Let us first show by induction that $p\in Int(C_{p})$. If $n=1$ the result is obvious. Let $n>1$. If $p\in Int(Q)$ for some cube $Q$
then the result follows immediately.
If $p\in\partial Q$ for some cube $Q$ then $p$ belongs to at least one  $(n+1)$-face $F$ of $Q$. The cubulation ${\cal{C}}$ induces a cubulation
${\cal{C}}_{n+1}$ in the hyperplane which contains the face $F$. By induction hypothesis $p$ is in the interior (relative to the hyperplane) of the union of
the $(n+1)$-cubes of ${\cal{C}}_{n+1}$ which contain $p$. Each of these $(n+1)$-cubes is a face of exactly two $(n+2)$-cubes of ${\cal{C}}$ which contain $p$.
Then $p$ is in the interior of the union of these cubes. Therefore, by induction, $p\in Int(C_{p})$.
To prove 1 it is sufficient to observe that, by the preceding argument, $p\in Int(C_{p})$, since  $C_{p}\subset{\cal {Q}}_{P}$. \\

\hspace{-.67cm}{\it 2.} Let $P_{m}$ be the hyperplane parallel to $P$ at distance $m$ of $P$ which intersects $Q$. Then $P_{m}$ is a support plane of
$Q$ and therefore there exists
a linear functional $\alpha:{\mathbb R}^{n+2}\to \mathbb R$ such that  $P_{m}=\alpha^{-1}(\{m\})$ and
$\alpha(p)\geq{m}$ for all $p\in{Q}$.
Since $Q$ is a convex polytope, it follows from standard facts of the geometry of convex sets and linear programming that $Q\cap P_{m}$ is a face of $Q$
(cubic simplex) of some dimension since it is the set where the linear function $\alpha$, restricted to $Q$, achieves its minimum.\\

\hspace{-.67cm}{\it 3.} We have that $\partial {\cal {Q}}_{P}\subset\underset{Q\subset\cal{Q}_{P}}\cup{\partial Q}$. Therefore
$\partial {\cal {Q}}_{P}$ is contained in a union of $(n+1)$-faces. Each face $F$ of a cube in $\cal C$ is a face of exactly two cubes of $\cal C$.
Furthermore $\partial {\cal {Q}}_{P}$ consists of
faces $F$ of $(n+2)$-cubes in $\cal {Q}_{P}$ with the property that $F$ is also the face of a cube not belonging to  ${\cal {Q}}_{P}$.
This is true since such type of faces belong to $\partial {\cal {Q}}_{P}$ and every point in $\partial {\cal {Q}}_{P}$ is contained in one of
those faces by the proof of 1 in proposition 2.2.
$\blacksquare$

\begin{lem}
Let $P$ be a hyperplane in $\mathbb{R}^{n+2}$. Let $\cal{C}$ be the canonical cubulation of $\mathbb{R}^{n+2}$ and
${\cal {Q}}_{P}=\cup\{Q\in {\cal{C}}:Q\cap P\neq\emptyset\}$.
Let $k$ be a point in $P$ and
$L_{k}$ the normal line to $P$ at $k$. Then $J_{k}=L_{k}\cap {\cal {Q}}_{P}$ is connected. By proposition 2.2,
$P$ is contained in the interior of ${\cal {Q}}_{P}$ and therefore $J_k$ is a non-trivial compact interval.
\end{lem}

{\it Proof.} $P$ divides $\mathbb{R}^{n+2}$ in two open connected components $H^+$ and $H^-$. Let us consider one of these components,
for instance $H^+$. Let $Q\in\cal{C}$ and suppose that $Q\subset{H^+}$. Then $Q\cap{P}=\emptyset$ and therefore $d(Q,P)=m>0$. By proposition 2.2 the set
${\cal{M}}:=\{q\in{Q}\,\,\,|\,\,\,d(q,P)=m\}$ is a cubic simplex, of some dimension contained in the boundary of $Q$. For each $(n+1)$-hyperface $F$
of $Q$ which intersects ${\cal{M}}$, there is a hyperplane $P_{F}$ which supports $F$. Only one of the closed halfspaces determined by
$P_{F}$ contains $Q$. We define $S(Q)$ as the intersection of all such closed halfspaces (see Figure 1).

\begin{figure}[tbh]
\centerline{\epsfxsize=2.2in \epsfbox{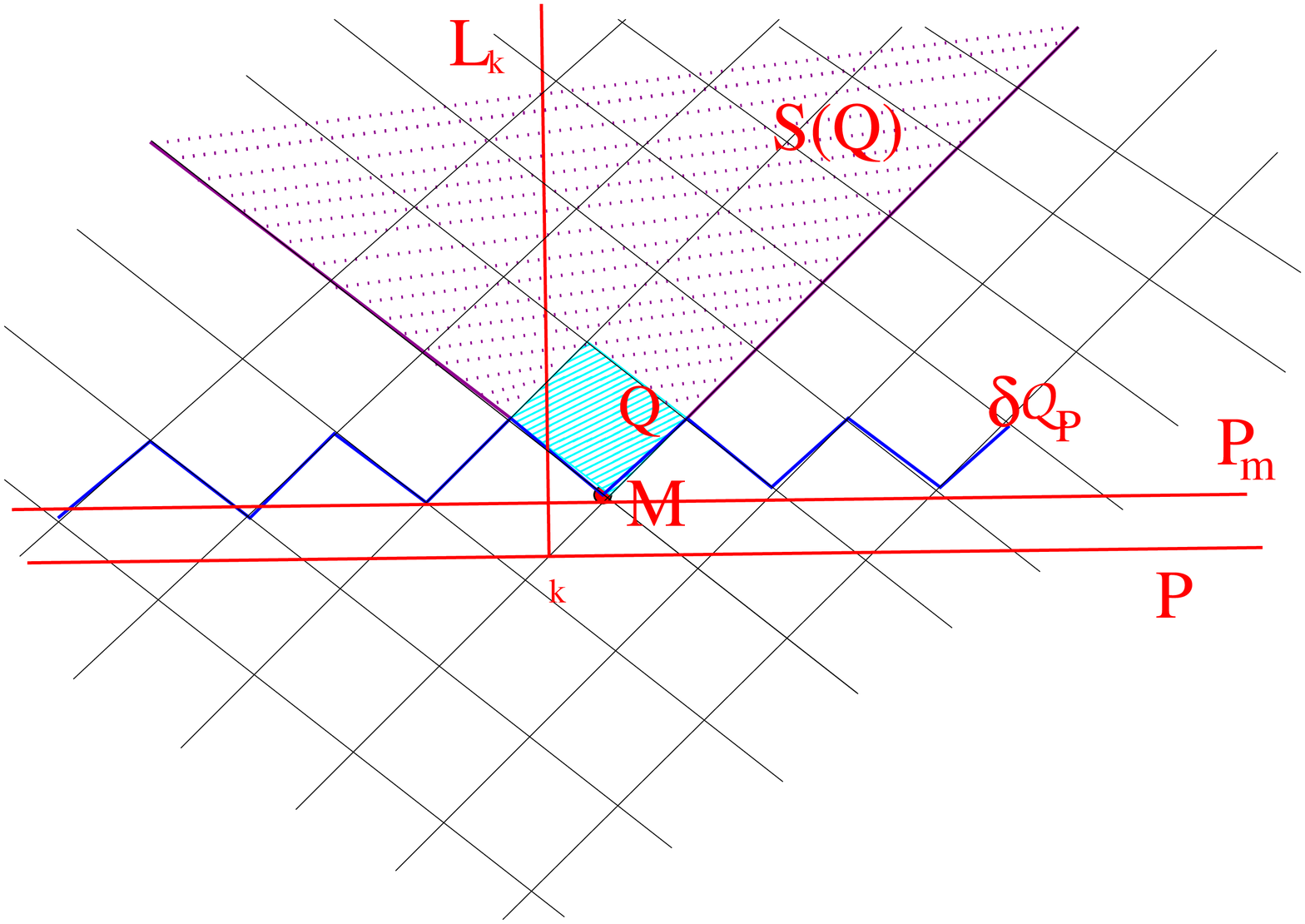}}
\caption{\sl The convex set $S(Q)$.}
\label{F1}
\end{figure}

The set $S(Q)$ is a convex unbounded set which contains $Q$ and it is a union of cubes of $\cal C$. We will show that $P_m$ is a support hyperplane of $S(Q)$
and $P_m\cap{S(Q)}={\cal{M}}$.\\

First we will prove:\\

\noindent {\bf Claim 1.} $S(Q)\cap P=\emptyset$. \\

{\it Proof of Claim 1}.  If $\cal{M}$ is a $(n+1)$-dimensional face
of $Q$, the result follows easily. If $\mbox{dim}({\cal{M}})\,<\,n+1$, then $\cal{M}$  is contained in the intersection of two hyperplanes
$P_{F_{1}}$ and $P_{F_{2}}$. One has that $P_{F_{1}}\cap P_{F_{2}}=P_{F_{1}}\cap P_{m}$ is of dimension $n$, and therefore divides $P_{F_{1}}$ into
two components, only one of which contains $Q\cap P_{F_{1}}$. Therefore $S(Q)$ lies in the halfspace determined by $P_{m}$ which
contains $Q$ and this halfspace does not contain $P$.  This proves claim 1.\\

To continue with the proof let us describe $P$, $P_m$ and $S(Q)$ in terms of linear equalities and linear inequalities.\\

Let us consider the set of hyperplanes $P_{F_{i}} \,\,(i=1,\ldots,k)$ defining $S(Q)$ which are not perpendicular to $P$.
Then the intersection of the halfspaces corresponding to these hyperplanes and containing $Q$, is defined by the set of points
$x\in{\mathbb R}^{n+2}$  which satisfy the following set of inequalities:

\setcounter{equation}{0}
\begin{equation}
\begin{array}{lll}
\langle x,{\text {\bf n}}_1\rangle &\geq&{a_1}\\
\langle x,{\text {\bf n}}_2\rangle &\geq&{a_2}\\
       &\vdots&\\
\langle x,{\text {\bf n}}_{k}\rangle &\geq&{a_{k}},
\end{array}
\end{equation}

\noindent where $\langle\cdot,\cdot\rangle$ is the standard inner product and $\langle {\text {\bf n}}_i,{\text {\bf n}}_j\rangle =\delta_{ij}$ and $a_i\in\mathbb R$.
If none of the hyperplanes defining $S(Q)$ is orthogonal to $P$ then $k=n+2$.\\

\noindent{\bf Definition:} The vector ${\text {\bf n}}_i$ is called the {\it exterior normal vector} to the face ${F_{i}}\subset Q$ contained in the hyperplane $P_{F_{i}}$.\\

If a hyperplane $P_{F_{j}}$ is perpendicular to $P$ the hyperplane corresponding to the opposite face to $F_{j}$ is also
perpendicular to $P$. Let $H_{F_{j}}$ be the closed halfspace determined by  $P_{F_{j}}$ which contains $Q$. \\

\noindent {\bf Obvious Remark.} If $L$ is a line which is perpendicular to $P$ then either $L\cap{H_{F_{j}}}=\emptyset$ or $L\cap{H_{F_{j}}}=L$.\\

Using a translation, if necessary, we can assume without loss of generality that the hyperplane $P$ is given by the linear equation:

\setcounter{equation}{1}
\begin{equation}
P=\{x\in{\mathbb R}^{n+2} \,\,\,\,|\,\,\langle x,{\text {\bf n}}\rangle =0\} \,\,\,\text{where}\,\, |{\text {\bf n}}|=1,
\end{equation}

\noindent where ${\text {\bf n}}$ is chosen in such a way that  $H^{+}=\{x\in{\mathbb R}^{n+2}\,\,\,|\,\,\langle x,{\text {\bf n}}\rangle >0\}$ (see Figure 2).

\begin{figure}[tbh]
\centerline{\epsfxsize=1.7in \epsfbox{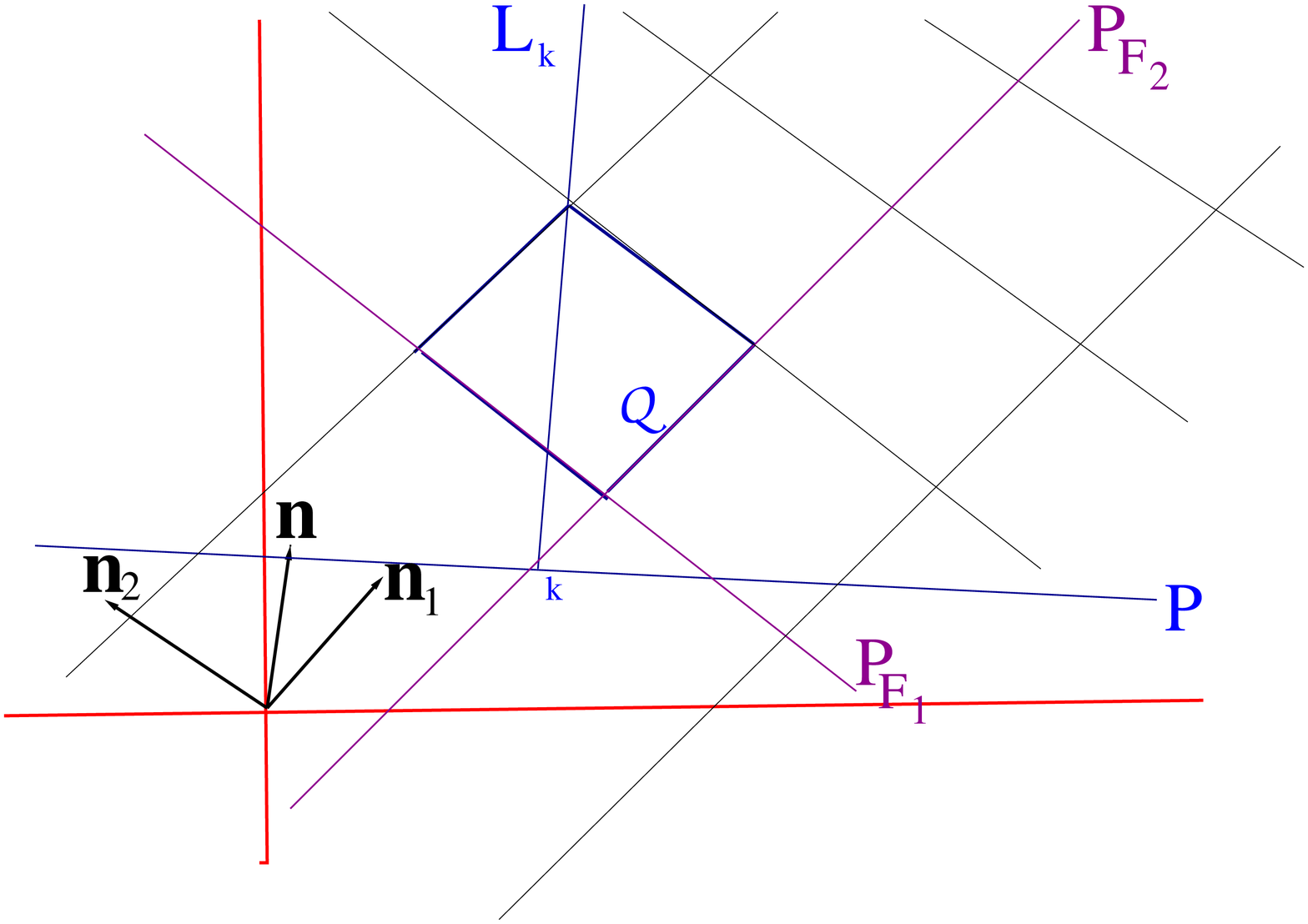}}
\caption{\sl Description of the set $S(Q)$ by inequalities.}
\label{F2}
\end{figure}

The hyperplane $P_{m}$ divides
$\mathbb{R}^{n+2}$ into two closed halfspaces, one of which contains $Q$ and the other contains $P$.\\

Then the hyperplane $P_{m}$ is given by the equation:
\setcounter{equation}{2}
\begin{equation}
P_{m}=\{x\in{\mathbb R}^{n+2} \,\,\,\,|\,\,\langle x,{\text {\bf n}}\rangle =m\},
\end{equation}
\noindent and by hypothesis we have that $\langle {\text {\bf n}},{\text {\bf n}}_i\rangle >0,\,\,i=1,\ldots,k$.\\

For each $p\in\Bbb R^{n+2}$, let $L_{p}$ be the normal line to $P$ which contains $p$.
Then the line $L_{p}$ can be parametrized by the function $\gamma:\mathbb R\to {\mathbb R}^{n+2}$ defined as follows:
\setcounter{equation}{3}
\begin{equation}
t\mapsto p+t{\text {\bf n}}
\end{equation}
{\bf Claim 2.}  $L_{p}\cap{S(Q)}$ is either empty
(it happens only if $Q$ has a hyperface, of dimension $n+1$, which is contained in a hyperplane perpendicular to $P$) or $L_{p}\cap{S(Q)}$ is a ray
{\it i.e.,} a noncompact closed interval contained in $L_{p}$. \\

{\it Proof of Claim 2.} By the obvious remark above we need to show that if $L_{p}\cap{S(Q)}\neq\emptyset$ then
$L_{p}\cap{S(Q)}=\{\gamma(t)\,\,|\,\,t\geq{c}\}$ where $c\geq{m}>0$. But this fact is immediate since, by (1) and (4),
it follows that if  $\gamma(t)\in L_{p}\cap{S(Q)}$ then since  $\langle {\text {\bf n}},{\text {\bf n}}_i\rangle >0,\,\,i=1,\ldots,k$ we have that
$t\geq\frac{a_j-\langle p,{\text {\bf n}}_j\rangle}{\langle{{\text {\bf n}}},{\text {\bf n}}_j\rangle}=c_j$.

We can choose $c=\max\{c_j\}$. Since
$S(Q)$ is contained in the halfspace determined by $P_m$ that contains $Q$ we must have $c\geq{m}$. This finishes the proof of claim 2.\\

Claim 2 implies that for each $q\in{Q}$, $L_q\cap{S(Q)}$ is a ray. Let $T(Q)$ be the union of all such rays. By the above we have the following\\

{\bf Corollary of Claim 2}. If $Q'\in\cal{C}$ is a cube which intersects $Int(T(Q))$ then $Q'\cap{P}=\emptyset$.\\

Now we are able to finish the proof of lemma 2.3. For $k\in{P}$, let $J_{k}^{+}=L_{k}\cap H^+$ with the order induced by the distance to $P$.
Let $x$ be the first point on $J_{k}^{+}$
such that $x\in J_{k}^{+}\cap\partial{\cal {Q}}_{P}$.
By proposition 2.2 the boundary  $\partial{\cal {Q}}_{P}$
is the union of $(n+1)$-dimensional cubes, each of which is a hyperface of exactly two hypercubes in $\cal{C}$: one
intersecting $P$ and one completely contained in $H^+$. Let $C$ be the cube containing $x$ which does not intersect $P$. Then,
by the above corollary, for all $z>x$, all cubes containing $z$ do not intersect $P$. Hence, $x$
is the only point on $J_{k}^{+}$ which lies on $\partial{\cal {Q}}_{P}$. Therefore $J_{k}^{+}\cap {\cal {Q}}_{P}$ is connected. $\blacksquare$\\

Lemma 2.3 implies the following corollary:

\begin{coro}
The lengths of $J_k\subset{L_k}$  and $J_{k}^{+}\subset{L_k}$ vary continuously with $k$. These intervals have a natural order which
provides them with a continuous orientation.  We denote $J_k$  and $J_{k}^{+}$  by $[a(k),b(k)]$ and  $[k,b(k)]$, respectively. Furthermore, the functions
$\psi^{+}(k):=b(k)$ and $\psi^{-}(k):={a(k)}$ are homeomorphisms from $P$ to $\partial {\cal {Q}}_{P}\cap{H^+}$ and $\partial {\cal {Q}}_{P}\cap{H^-}$,
 respectively.

 As a consequence ${\cal {Q}}_{P}$ is a closed tubular neighborhood (or bi-collar neighborhood) of $P$.
\end{coro}

{\it Proof.} By lemma 2.3, $J_{k}^{+}\cap \partial{\cal {Q}}_{P}$ is connected. We will prove that $J_{k}^{+}\cap \partial{\cal {Q}}_{P}$ is in fact
one point. Since the cubes in the cubulation  $\cal{C}$ are right-angled, if $F$ is a face of a cube orthogonal to $P$ then $P$ intersects $F$. This
means that $F$ can not be contained in the boundary of ${\cal {Q}}_{P}$. Therefore $J_{k}^{+}$, which is orthogonal to $P$, must intersect each
$(n+1)$-cube in $\partial{\cal {Q}}_{P}$ transversally.\\

We have proved that $J_{k}$ is connected and in particular, $J_{k}\cap \partial {\cal {Q}}_{P}$ consists of two points
$a(k)$ and $b(k)$. Let $b(k)$ be the corresponding point in $H^{+}$. We define the map
$\psi^{+}:P\rightarrow \partial {\cal {Q}}_{P}\cap{H^+}$ such that
 $x\in P$ is sent to $b(x)$. Without loss of generality we can assume that
$$
P=\{x=(x_1,\ldots,x_{n+1},x_{n+2})\in {\mathbb R}^{n+2}\,\,\,|\,\,\,x_{n+2}=0\}.
$$
 Thus can identify $P$ with ${\mathbb R}^{n+1}$. In terms of this identification we have $\psi^{+}(k)=(k,h^+(k)), \,\,\,k\in{\mathbb R}^{n+1}$,
and $h^+(k)=d(k,\psi^+(k))=d(k,b(k))$.
Consequently $H^{+}\cap \partial {\cal {Q}}_{P}$ can be thought of as the graph of the function $h^{+}$ and since this graph is closed,
it follows that $h^+$ is continuous.  This implies that the length of $J^{+}_{k}$ varies continuously. Hence $\psi^{+}$ is a homeomorphism.
Obviously all the preceding arguments remain valid if we had used  $H^{-}$, $\psi^-$ and  $J_-$. $\blacksquare$

\begin{coro} Since ${\cal {Q}}_{P}$ is a bi-collared neighborhood of $P$ we obtain that $P$ can be continuously isotoped onto any
of the two connected components of the boundary $\partial {\cal {Q}}_{P}$. This can be done by an ambient isotopy of $\Bbb R^{n+2}$
\cite{cernavskii}, \cite{rushing}.
\end{coro}

We will use the following lemma to prove our main theorem.

\begin{lem}
 Let $P_{1}$ and $P_{2}$ be two orthogonal codimension one
hyperplanes in $\mathbb{R}^{n+2}$. Then $P_{1}\cap P_{2}$ can be cubulated;
more precisely, there is an ambient isotopy which takes $P_{1}\cap P_{2}$ onto a cubic complex,
contained in the scaffolding $\cal S$ of the canonical cubulation $\mathcal{C}$ of $\mathbb{R}^{n+2}$.
\end{lem}

{\it Proof.} Let $\mathcal{Q}_{1}$ be the cubic complex formed by the
union of elements of $\mathcal{C}$ intersecting $P_{1}$. By the
above corollary, $\mathcal{Q}_{1}$ is a tubular neighborhood (or bi-collar in this case)
of $P_{1}$. Let $E$ be a connected component of the boundary of $\mathcal{Q}_{1}$. Recall the isotopy $\psi=\psi^{+}$ constructed in corollary 2.4 from $P_{1}$
to $E$, which assigns to each point $x$ in $P_{1}$
the unique point on $E$ lying on the line normal to $P_{1}$ at $x$. The
inverse of $\psi $ in this case is in fact the canonical projection $\pi_{1}$
of $\mathbb{R}^{n+2}$ onto $P_{1}$.\\

Since $P_{1}$ and $P_{2}$ are orthogonal, if $x$ is a point in $P_{1}\cap
P_{2}$ then the normal line at $x$ to $P_{1}$ will lie in $P_{2}$.
Therefore, $\psi $ restricted to $P_{1}\cap P_{2}$ is a homeomorphism which can be extended to an ambient isotopy (\cite{cernavskii})
between $P_{1}\cap P_{2}$ and $E\cap P_{2}$.
Now by proposition 2.2 $E$ is a cubic complex of dimension $n+1$. We will repeat the idea of
lemma 2.3 for $P_{2}\cap
E\subset E$. Take the union $\mathcal{B}$ of all $(n+1)$-dimensional cubes
in $E$ that intersect $P_{2}\cap E$. We will prove that $\mathcal{B}$ is a
bi-collar of $P_{2}\cap E$ in $E$ and therefore one of its boundary components is ambient isotopic to $P_1\cap{P_2}$.
This boundary component is obviously contained in the $n$-skeleton. This would prove the lemma.\\

The intersection $P_{1}\cap P_{2}$ is an $n$-dimensional hyperplane, and
hence of codimension one in $P_{1}$. Consider the foliation of $P_{1}$ given
by lines orthogonal to $P_{1}\cap P_{2}$. Notice that the image of each of
these lines under $\psi $ is a polygonal curve $\gamma$ in $E$, intersecting $E\cap
P_{2}$ in one unique point. \\

As mentioned above, $E$ is cubulated as a union of $(n+1)$-faces of cubes of $\cal C$. The projection $\pi_{1}$ of this $(n+1)$-dimensional cubulation gives a decomposition
${\cal{P}}$ of $P_{1}$ into parallelepipeds. Observe that, in general, this decomposition is not a cubulation. We will prove
that a result analogous to the lemma 2.3 holds for ${\cal{P}}$. By this we mean the following: Let $Q$ be a parallelepiped in
${\cal{P}}$ which does not intersect  $P_{1}\cap P_{2}$, $x$ be a point in $Q$ and $l_{x}$ be the set of points $y$ on the line in
$P_1$ normal to
$P_{1}\cap P_{2}$ which passes through $x$ such that $x<y$, with the order induced by the distance to $P_2$.  We will prove that any $Q'\in {\cal{P}}$
which intersects $l_{x}$ does not intersect $P_{1}\cap P_{2}$.\\

The parallelepiped $Q$ is the image under $\pi_{1}$ of an $(n+1)$-cube in $E$. This $(n+1)$-cube  is a face of two $(n+2)$-cubes in
$\mathbb{R}^{n+2}$, at least one of which does not intersect $P_{2}$. Let $C$ be the $(n+2)$-cube not intersecting $P_{2}$. Then,
the set $S(C)$ defined in the proof of the lemma 2.3 does not intersect $P_{2}$. Thus $\pi_{1}(S(C))$ does not intersect
$P_{1}\cap P_{2}$. Moreover, $l_{x}$ is contained in $\pi_{1}(S(C))$ and,
any parallelepiped in ${\cal{P}}$ which intersects the interior of $\pi_{1}(S(C))$ is contained in $\pi_{1}(S(C))$, and
therefore does not intersect $P_{1}\cap P_{2}$. Hence, no cube in $\psi(\pi_{1}(S(C)))$ intersects $E\cap P_{2}$. That is, no
$(n+1)$-cube in $E$ which intersects the polygonal ray $\psi(l_{x})$ intersects $E\cap P_{2}$.\\

Proceeding as in corollaries 2.4 and 2.5, choosing one of these points for each $x\in E\cap
P_{2}$ yields a continuous function from $E\cap P_{2}$ to the boundary of $\mathcal{B}$.
Hence, $\mathcal{B}$ is a bi-collar of $E\cap P_{2}$ in $E$ and  $E\cap P_{2}$ can be deformed by an ambient isotopy into
one boundary component of $\mathcal{B}$. Observe that this isotopic copy of $E\cap P_{2}$ is contained in the $n$-dimensional skeleton of the cubulation
$\mathcal{C}$. $\blacksquare $

\begin{rem} In lemmas 2.3 and 2.6 we only need to consider the subset of those cubes of $\mathcal C$ whose
distance to the corresponding hyperplanes $P$, $P_1$ and $P_2$ is sufficiently small, for instance, less or equal
than $4\sqrt{n+2}$, {\it i.e.,} the cubes between two parallel hyperplanes at distance  $8\sqrt{n+2}$. The number $\sqrt{n+2}$
appears because it is the diameter of the unit cubes in $\Bbb R^{n+2}$.
\end{rem}

A modification of the methods of the preceding results for hyperplanes can extend lemma 2.3 to the more general case given by following theorem:

\begin{theo}
Let  $M^{n+1}\subset\mathbb{R}^{n+2}$ be a smooth, compact and closed manifold. Let $V(M)$ be a closed tubular neighborhood of $M$.
We can assume that $V$ is the union of linear segments of equal length $c>0$, centred at
points of $M$, and orthogonal to $M$. Let $\cal{C}$ be the canonical cubulation
of $\mathbb{R}^{n+2}$. For $m\in\Bbb N$, let $\mathcal {C}_{m}$ denote the corresponding subcubulation.

Let ${\cal {Q}}_{M}=\cup\{Q\in {\mathcal{C}_m}:Q\cap M\neq\emptyset\}$.
Then, we can choose $m$ big enough such that ${\cal {Q}}_{M}$ is a closed tubular neighborhood of $M$ and $M$ can be deformed by an
ambient isotopy onto any of the two boundary components of this tubular neighborhood.
\end{theo}

{\it Proof.} The idea of the proof is to ``blow up''  the manifold and its tubular neighborhood by a homothetic transformation so that
inside balls of large (but fixed) radius the manifold is almost flat and the normal segments are almost parallel. Then, locally,
the manifold is approximately a hyperplane and we can apply a modification of the methods of
the previous lemmas and then we rescale back to the original size by the inverse homothetic transformation. The modification consists in replacing
the normal segments to a hyperplane by segments of the flowlines of a nonsingular vector field which is very close to these normal 
segments in the tubular neighborhood of $M$.\\

More precisely, let $\frak{h}_m:{\Bbb R}^{n+2}\to{\Bbb R}^{n+2}$ be the homothetic transformation $\frak{h}_m(x)=mx$
where $m$ is a positive integer. Let $M_m=\frak{h}_m(M)$ and $V_m=\frak{h}_m(V)$.  The homothetic transformation $\frak{h}_m$ is isotopic to the identity,
hence $M_m$ is isotopic to $M$ and $V_m$ is isotopic to $V$. Given
$\epsilon>0$ we can choose $m$ large enough such that at each point of $M_m$ its sectional curvature is less than $\epsilon$.
This  means we can choose $m$ large enough such that for every point $p\in{M_m}$ the pair $(B(p),B(p)\cap{M_m})$ is
$C^\infty$-close to the pair $(B(p),B(p)\cap{T_pM_m)}$, where $B(p)$ denotes the closed ball centred at $p$ of radius sufficiently large, for instance of radius
$10m\sqrt{n+2}$, and $T_pM_m$ is the tangent space of $M_m$ at $p$. \\

Let  $x\mapsto {\text {\bf n}}_{x}$, $x\in M_{m}$ be the unit normal vector field of $M_{m}$ with respect to an orientation of $M_{m}$.  It follows from standard
facts of differential geometry (see \cite{milnor}) that the normal map $x\mapsto{x+t{\text {\bf n}}_{x}}$, for $x\in B(p)\cap{M_m}$
has no focal points in  $V_m$. This  implies that the map $\mu_{p}:(B(p)\cap{M_m})\times [-1,1]\rightarrow V_{m}$
given by $(x,t)\mapsto{x+t{\text {\bf n}}_{p}}$
is a diffeomorphism onto its image. Furthermore, there exists a constant $\delta_p>0$ such that if $\text{\bf v}_p$ is a unit vector
satisfying $||{\text{\bf v}}_p-{\text{\bf n}}_{p}||<\delta_p$ then the map $x+t{\text{\bf v}}_p$ is still a diffeomorphism from
$(B(p)\cap{M_m})\times [-1,1]$ onto its image and $\langle{{\text{\bf v}}_p,{\text{\bf n}}}\rangle >0$ for any normal vector to $M_m$ at a point $x\in B(p)\cap{M_m}$. \\

The map $\varphi:M_{m}\times[-mc/2,mc/2]\to{V_{m}}\subset{\mathbb R}^{n+2}$
given by $\varphi(x,t)=x+t{\text {\bf n}}_{x}$ is a parametrization of $V_{m}$.
Let $V_{1/2}=\varphi(M_{m}\times[-mc/4,mc/4])$ be the smaller neighborhood of width $mc/2$. The boundary of $V_{1/2}$ has two connected components
$$
\partial V_{1/2}^{-}=\varphi(M\times\{-mc/4\})
$$
$$
\partial V_{1/2}^{+}=\varphi(M\times\{mc/4\}).
$$

Let ${\cal {Q}}_{M_{m}}=\cup\{Q\in {\mathcal{C}}:Q\cap M_{m}\neq\emptyset\}$. Its boundary has also two connected components 
$\partial{\cal {Q}}_{M_m}^{+}:=\partial{\cal {Q}}_{M_m}\cap V_{1/2}^{+}$ and
$\partial{\cal {Q}}_{M_m}^{-}:=\partial{\cal {Q}}_{M_m}\cap V_{1/2}^{-}$.\\

\noindent{\bf Claim.} There exists a finite family of diffeomorphisms
$$\{\psi_i:{\mathbb B}^{n+1}\times[-1,1]\to{V_{m}}\}_{i=1}^{k}$$
\noindent such that if we set $U_i:=\psi_i({\mathbb B}^{n+1}\times[-1,1])$ then
\begin{enumerate}

\item $V_{1/2}\subset\cup_{i=1}^{k}\,\,Int(U_i)$

\item $\psi_i({\mathbb B}^{n+1}\times\{0\})\subset{M_{m}}$

\item  For each $i$ the curves  $t\mapsto\psi_i(y,t)$ are parallel to the unit vector ${\text{{\bf v}}}_i$

\item  For every cube $Q$ of $\cal C$ such that $Q\cap{M}_m=\emptyset$, $Q\cap\partial{\cal {Q}}_{M_m}^{+}\neq\emptyset$
and $Q\subset{Int(U_i)}$ we have that $\langle {\text{{\bf v}}}_i, w  \rangle>0$  
for every $w$ exterior normal to a face of $Q$
which is contained in $\partial{\cal {Q}}_{M_m}^{+}$ (see definition in proof of lemma 2.3).

\item We define the set $S(Q)$ as in proof of lemma 2.3, where the set $\cal{M}$ is replaced
by $Q\cap\partial{\cal {Q}}_{M_m}^{+}$ (see also equation 1). We consider the convex set $S_{V_{m}}(Q)=S(Q)\cap V_{m}$. Then there exists $t_{0}\in (-1,1)$ such that
$\psi_i(y,t)\in S_{V_{m}}(Q)$, $t\geq t_{0}$.
\end{enumerate}

{\it Proof of Claim.} Since $M_m$ is compact, there exist a finite open subcovering
$\{A_{1},A_{2},\,\ldots\,,A_{k}\}$ of the open covering $\{Int(B(p))\cap{M_m}\}_{p\in{M_m}}$, where $A_{i}:=Int(B(p_{i}))\cap{M_m}$ and
vectors ${\text{\bf v}}_{p_{1}},\ldots, {\text{\bf v}}_{p_{k}}$ such that
$||{\text{\bf v}}_{p_i}-{\text{\bf n}}_{p_i}||<\delta_{p_i}$ and the maps  $\tilde{\mu}_{p_{i}}:A_{i}\times [-1,1]\rightarrow V_{m}$
given by $(x,t)\mapsto x+t{\text {\bf v}}_{i}$ have the property that
$V_{1/2}\subset\cup_1^k\,\,\tilde\mu(A_i\times[-1,1])$. Furthermore,  ${\text{\bf v}}_{p_{1}},\ldots, {\text{\bf v}}_{p_{k}}$  
can be chosen such that $\langle {\text{{\bf v}}}_i, w  \rangle>0$  for every $w$ exterior normal to a face of $Q$
which is contained in $\partial{\cal {Q}}_{M_m}^{+}$ and intersects $\mu_{p_{i}}(A_{i}\times [-1,1])$ (see Figure 3).
 Notice that this can be done
because $Q$ is a cube, so its faces meet at right angles.
Using the fact that $A_{i}$ is homeomorphic to $Int({\mathbb B}^{n+1})$, the claim follows.

\begin{figure}[tbh]
\centerline{\epsfxsize=2in \epsfbox{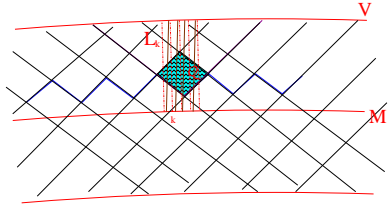}}
\caption{\sl Brown lines  $t\mapsto\psi_i(y,t)$ parallel to the same unit vector ${\text{{\bf v}}}_i$.}
\label{F3}
\end{figure}

\noindent{\bf Remark:} The above claim can be modified to prove a similar result for  $\partial{\cal {Q}}_{M_m}^{-}$.\\

For each $i=1, \ldots,k$, let $f_i:{\mathbb R}^{n+2}\to[0,1]$  be a smooth function such that
$f_i(x)>0$ if $x\in Int(U_i)$ and $f_i(x)=0$ if $x\notin\,Int(U_i)$.
Let ${\frak v}_i(z)={\text{{\bf v}}}_i$ be the global constant vector field equal to ${\text{{\bf v}}}_i$ at every point
 $z\in{\mathbb R}^{n+2}$. Let $\frak{V}(z)=\sum_{i=1}^k\,\,f_i(z){\text{{\bf v}}}_i$.\\

Then $\frak{V}(z)$ is nonzero in $V_{1/2}$ and it has compact support inside $V$, hence it defines a global flow $\eta_t:{\mathbb R}^{n+2}\to{\mathbb R}^{n+2}$.\\

\noindent{\bf Claim.}  Let $\partial V_{1/2}^{-}$ and $\partial V_{1/2}^{+}$ be the two boundary components of $V_{1/2}$.
Then
\begin{enumerate}
 \item  If $z\in \partial V_{1/2}^{-}$ there exits $t_{z}>0$ such that $\eta_{t_{z}}(z)\in\partial V_{1/2}^{+}$.

\item  The  flowlines of $\{\eta_t\}_{t\in\mathbb R}$ are transversal to $\varphi(M\times\{t\}),\,\,\,t\in[-c/4,c/4]$.

\item The  flowlines of $\{\eta_t\}_{t\in\mathbb R}$ meet both
$\partial V_{1/2}^{-}$ and $\partial V_{1/2}^{+}$ in one point and furthermore the flowlines meet both
$\partial{\cal {Q}}_{M_m}^{+}:=\partial{\cal {Q}}_{M_m}\cap V_{1/2}^{+}$ and
$\partial{\cal {Q}}_{M_m}^{-}:=\partial{\cal {Q}}_{M_m}\cap V_{1/2}^{-}$ in one point.
\end{enumerate}

{\it Proof of Claim.} First we remark that the flow is defined for all $t\in \mathbb R$ because the vector field  $\frak{V}$
has compact support. Every point of $V_m$ can be written uniquely in the form $x+t{\text{{\bf n}}}_x,\,\,x\in{M_m},\,\,t\in[-cm/2,cm/2]$.  
Let $\frak p:V_m\to[-cm/2,cm/2]$ be the map
$x+t{\text{{\bf n}}}_x\mapsto{t}$. By the properties of $\frak V$, we have that for fixed $y$ the function
$t\mapsto\frak p\circ\eta_t(y)$
is strictly increasing for $t\in[-cm/4,cm/4]$ since the derivative of the function is positive. This implies that the flowlines meet both 
$\partial V_{1/2}^{+}$ and $\partial V_{1/2}^{-}$ in one point.
Furthermore, if $Q$ is a cube which meets  $\partial{\cal {Q}}_{M_m}$ but does not meet $M_m$ we have that $S(Q)$, defined above, is 
positively invariant under the flow, in fact $\eta_t(S(Q))\subset {Int(S(Q))}$. Therefore
the flowlines meet both $\partial{\cal {Q}}_{M_m}^{+}$ and
$\partial{\cal {Q}}_{M_m}^{-}$ in exactly one point. This proves the claim.\\

We can rescale the flow by multiplying the vector field $\frak V$ by a positive smooth function to obtain a new flow
$\{{{\hat{\eta}}_t}\}_{t\in\mathbb R}$ such that
for every  $z\in{M_m}$ we have $\hat\eta_{-1}(z)\in\partial{\cal {Q}}_{M_m}^{-}$  and
$\hat\eta_{1}(z)\in\partial{\cal {Q}}_{M_m}^{+}$.\\

Then the map $\Phi:M\times[-1,1]\to{\mathbb R}^{n+2}$ defined by $\Phi(x,t)=\hat\eta_t(x)$ has the properties:
\begin{enumerate}
\item $\Phi(M\times[-1,1])={\cal {Q}}_{M_m}$ and therefore ${\cal {Q}}_{M_m}$ is a closed tubular neighborhood of $M$.
\item $\Phi(M\times\{-1\})=\partial{\cal {Q}}_{M_m}^{-}$
\item $\Phi(M\times\{1\})=\partial{\cal {Q}}_{M_m}^{+}$
\end{enumerate}

We have concluded that the set of cubes of the
{\it canonical} cubulation of $\mathbb R^{n+2}$ that touch $M_m$ is a closed tubular neighborhood of $M_m$ and $M_m$ can be
deformed by an ambient isotopy  onto any of the two boundary components of this tubular neighborhood.
To finish the proof we now rescale our construction back to its original size using the inverse homothetic transformation $\frak{h}_{\frac1m}$.
This transformation transforms $\mathcal C$ onto the subcubulation ${\mathcal C}_m$.  $\blacksquare$\\

Let  $M,\,\,N\subset\Bbb R^{n+2}$, $N\subset{M}$,
be compact, closed and smooth submanifolds of  $\Bbb R^{n+2}$ such that $\mbox{dimension}\,(M)=n+1$ and $\mbox{dimension}\,(N)=n$.
Since $M$ is codimension one in $\Bbb R^{n+2}$ it is oriented
and we will assume that $N$ has a trivial normal bundle in $M$ ({\it i.e.,} $N$ is a two-sided hypersurface of $M$). Then under
these hypotheses we have the following theorem for pairs $(M,N)$ of smoothly embedded compact submanifolds of $\Bbb R^{n+2}$:

\begin{theo} There exists an ambient isotopy of $\Bbb R^{n+2}$ which takes $M$ into the $(n+1)$-skeleton of the canonical
cubulation $\cal C$ of $\Bbb R^{n+2}$ and $N$ into the $n$-skeleton of $\cal C$. In particular, $N$ can be deformed by
an ambient isotopy into a cubical manifold contained in the canonical scaffolding of  $\Bbb R^{n+2}$.
\end{theo}

{\it Proof.} The proof is very similar to that of the previous theorem and lemma 2.6. As before, given $\epsilon>0$ there exists $m\in\Bbb N$,
large enough, such that
if we consider the homothetic transformation  $\frak{h}_m(x)=mx$  then the sectional curvatures of {\it both} $\frak{h}_m(M):=M_m$ and $\frak{h}_m(N):=N_m$ are
less than $\epsilon$. More precisely, for every $p\in{N_m}$ the triple
$(B(p),B(p)\cap{M_m}, B(p)\cap{N_m})$ is $C^\infty$-close to the triple $(B(p),B(p)\cap{T_pM_m}, B(p)\cap{T_pN_m})$,
where $B(p)$ denotes the closed ball centred at $p$ of radius sufficiently large, for instance of radius
$10m\sqrt{n+2}$, and $T_pM_m$ is the tangent space of $M_m$ at $p$ and $T_pN_m\subset{T_pM_m}$ is the tangent space to $N_m$ at $p$.\\

Let ${\mathcal Q}_{M_m}$, $V_m$, $V_{1/2}$, $\{{{\hat{\eta}}_t}\}_{t\in\mathbb R}$, and $\Phi$ be as in the previous theorem.  Let $f^+:M_m\to\partial{\cal Q}_{M_m}^+$ be 
given by $f^+(z)=\Phi(z,1)= \hat\eta_1(z)$. Let $\hat{N}_m=f^+(N_m)$. From above we know that
$\partial{\cal {Q}}_{M_m}^+$ is isotopic to $M_m$ and it is cubulated since it is a union of faces of cubes of $\cal C$.
We have that $N_m$ can be deformed to $\hat{N}_m$ by a global isotopy of ${\mathbb R}^{n+2}$. Since $\hat{N}_m$ is contained in 
$\partial{\cal {Q}}_{M_m}^+$ to prove the theorem we only need
to find an isotopy $\{h_t:{\mathbb R}^{n+2}\to{\mathbb R}^{n+2}\}_{t\in[0,1]}$ such that
$h_t(\partial{\cal {Q}}_{M_m}^+)\subset\partial{\cal {Q}}_{M_m}^+$,
$h_0=\text{Identity}$ and $h_1(\hat{N}_m)$ is contained in the $n$-skeleton of $\cal C$.\\

Let ${\cal Q}_{\hat{N}_m}$ be the union of all faces of $\partial{\cal {Q}}_{M_m}$ which intersect $\hat{N}_m$.
Let $\partial{\cal Q}_{\hat{N}_m}$ be the boundary of  ${\cal Q}_{\hat{N}_m}$ as a subset of $\partial{\cal Q}_{M_m}$. Then 
$\partial{\cal Q}_{\hat{N}_m}$ is a union of $n$-cubes.\\

\noindent{\bf Claim.}  ${\cal Q}_{\hat{N}_m}$ is a closed tubular neighborhood of  $\hat{N}_m$ in $\partial{\cal {Q}}_{M_m}$.\\

{\it Proof of Claim.} The proof is similar to the proof of lemma 2.6. Let  $A:N_m\times[-c,c]\to{M_m}$ be the parametrization of a tubular
neighborhood of $N_m$ in $M_m$ given by $A(x,s)=\gamma_x(s)$, where  $\gamma_x:\mathbb R\to{M_m}$ is the geodesic in $M_m$ 
(with respect to the induced metric of ${\mathbb R}^{n+2}$ in $M_m$) such that $\gamma_x(0)=x$ and 
$\gamma'_x(0)={\text{\bf w}}_x$ where ${\text{\bf w}}_x$
is a unit vector at $x$ tangent to $M_m$ and transversal to  $T_xN_m$. In other words we use a sort of Fermi coordinates for the tubular neighborhood. Given $d>0$ we can choose $m$  
large enough so that the curvature of both $M_m$ and $N_m$ is very small and the geodesics $\gamma_x(s)$ are very close (in the $C^2$-topology) 
to linear segments for $s\in[-d,d]$. We can take $c=d$. 
Let $y=f^+(x)\in \hat{N}_m$, $x\in N_m$ and let $\beta:[-c,c]\to \partial{\cal {Q}}_{M_m}^+$
the function $\beta(s)=f^+(\gamma_x(s))$ and $J_y=\{\beta(s)\,\,|\,\,s\in[0,c]\}$. Then just as in the proofs of lemma 2.6 and theorem 2.8 we can prove that we can choose the vector field $x\mapsto{\text{\bf w}}_x$ in such a way that
$J_y$ is homeomorphic to a non-trivial segment and intersects $\partial{\cal Q}_{\hat{N}_m}$ in exactly one point (see Figure 4). The curves $J_y$ are
 rectifiable and the function $y\mapsto\text{length}(J_y)$ is continuous. Hence, ${\cal Q}_{\hat{N}_m}$ is a bi-collar of $\hat{N}_m$ in 
$\partial{\cal {Q}}_{M_m}$. This proves the claim. 

\begin{figure}[tbh]
\centerline{\epsfxsize=2in \epsfbox{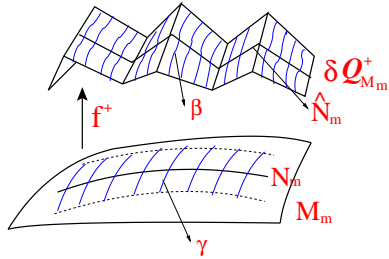}}
\caption{\sl The manifold $\hat{N}_m$ and the flow $\beta$.}
\label{F4}
\end{figure}

By the above,  $\hat{N}_m$ can be deformed  to one connected component of $\partial {\cal Q}_{\hat{N}_m}$ by an isotopy in $\partial{\cal {Q}}_{M_m}^+$. By standard
theorems (see  \cite{cernavskii}, \cite{rushing}), this
isotopy can be extended to a global isotopy 
$\{h_t:{\mathbb R}^{n+2}\to{\mathbb R}^{n+2}\}_{t\in[0,1]}$. 
Observe that this isotopic copy of $\hat{N}_m$ is contained in the $n$-dimensional skeleton of the cubulation
$\mathcal{C}$. To finish the proof we now rescale our construction back to its original size using the inverse
homothetic transformation $h_{\frac1m}$.
This transformation transforms $\cal C$ onto the subcubulation ${\cal C}_m$. $\blacksquare$

\section{Cubic knots}

In classical knot theory, a subset $K$ of a space $X$ is a {\it knot} if $K$ is homeomorphic to a sphere
$\mathbb{S}^{p}$. Two knots $K$, $K'$ are {\it equivalent} if there is a homeomorphism $h:X\rightarrow X$ such that $h(K)=K'$;
in other words $(X,K)\cong (X,K')$. However, a knot $K$ is sometimes defined to be an embedding
$K:\mathbb{S}^{p}\rightarrow\mathbb{S}^{n}$ or $K:\mathbb{S}^{p}\rightarrow\mathbb{R}^{n}$ (see \cite{mazur}, \cite{rolfsen}).
We are mostly interested in the codimension two smooth case $K:\mathbb{S}^{n}\rightarrow\mathbb{R}^{n+2}$.\\

In this section we will prove the main theorem of this paper.

\begin{theo} Let $\cal{C}$ be the canonical cubulation of $\mathbb{R}^{n+2}$. Let $K\subset \mathbb{R}^{n+2}$ be a smooth knot of dimension $n$.
 Given any closed tubular neighborhood $V(K)$ of $K$ there exists an ambient isotopy
$f_t:\mathbb{R}^{n+2}\to{\mathbb{R}^{n+2}}$ with support in $V(K)$ and $t\in[0,1]$ such that $f_0=Id$ and $f_1(K):=\bar{K}$ is contained
in the $n$-skeleton of the subcubulation ${\cal C}_m$ for some integer $m$. In fact $\bar{K}$ is contained in the boundary of the cubes of ${\cal C}_m$
which are contained in $V(K)$ and intersect $K$.
In particular, using the homothetic transformation $h_m(x)=mx$ we see that there exists a knot
$\hat{K}$ isotopic to $K$, which is contained in the scaffolding ($n$-skeleton) of the canonical cubulation  $\cal{C}$
of $\mathbb{R}^{n+2}$.
\end{theo}

{\it Proof}. Let $\Psi:\mathbb{S}^{n}\rightarrow \mathbb{R}^{n+2}$
be a smoothly embedded knotted $n$-sphere in $\mathbb{R}^{n+2}$. We denote $K=\Psi(\mathbb{S}^{n})$ and endow it with the
Riemannian metric induced by the standard Riemannian
metric of $\mathbb{R}^{n+2}$.\\

Let $V(K)$ be a closed tubular neighborhood of $K$. By a theorem of Whitney (\cite{whitney}) any embedded $\mathbb{S}^{2}$ in
$\mathbb{S}^{4}$ has trivial normal bundle. Since $H^{2}(\mathbb{S}^{n},\mathbb{Z})=0$ for $n>2$, any embedding of $\mathbb{S}^{n}$
in $\mathbb{R}^{n+2}$ has trivial normal bundle.  Hence $V(K)$ is
diffeomorphic to $K\times\mathbb{D}^{2}$ . Let $\phi: K\times\mathbb{D}^{2}\rightarrow V(K)$ be a diffeomorphism and  $p:V(K)\rightarrow K$ be the projection.
We can assume that the fibers of $p$ are Euclidean disks of radius
$\delta>0$. If $0<r<\delta$, the restriction of $\phi$ to $K\times \mathbb{D}_{r}^{2}$ (where $\mathbb{D}_{r}^{2}$
is the closed disk of radius $r$) is a closed
tubular neighborhood $W(K)$ of $K$. Then $K\subset W(K)\subset Int(V(K))$.
Take a section in $\partial W(K)$, \textit{i.e.} for a fixed $\theta _{0}\in \mathbb{S}_{r}^{1}$ consider
$\widetilde{K}=\phi (K\times \theta _{0})$. Observe that $\widetilde{K}$
is isotopic to the knot $K$. Now consider the pair $(\partial W(K),\widetilde{K})$. By the above $\widetilde{K}$ has normal
bundle in $\partial W(K)$ and applying theorem 2.9 yields
the result. $\blacksquare$

M. Boege. {\tt Instituto de Matem\'aticas, Unidad Cuernavaca}. Universidad Nacional Au\-t\'o\-no\-ma de M\'exico.
Av. Universidad s/n, Col. Lomas de Chamilpa. Cuernavaca, Morelos, M\'exico, 62209.

\noindent {\it E-mail address:} margaret@matcuer.unam.mx
\vskip .3cm
G. Hinojosa. {\tt Facultad de Ciencias}. Universidad Aut\'onoma del Estado de Morelos. Av. Universidad 1001, Col. Chamilpa.
Cuernavaca, Morelos, M\'exico, 62209.

\noindent {\it E-mail address:} gabriela@uaem.mx
\vskip .3cm
A. Verjovsky. {\tt Instituto de Matem\'aticas, Unidad Cuernavaca}. Universidad Nacional Au\-t\'o\-no\-ma de M\'exico.
Av. Universidad s/n, Col. Lomas de Chamilpa. Cuernavaca, Morelos, M\'exico, 62209.

\noindent {\it E-mail address:} alberto@matcuer.unam.mx


\begin{thebibliography}{99}

\bibitem{artin} E. Artin. \emph{Zur Isotopie zweidimensionalen
  Fl\"achen im $R_{4}$} Abh. Math. Sem. Univ. Hamburg
(1926), 174-177.
\bibitem{cairns} S. Cairns. \emph{A simple triangulation method for smooth manifolds}. Bull. Amer. Math. Soc. 67  pp. 389-390, 1961.
\bibitem{cernavskii} A. V. \v{C}ernavski\v{i}. \emph{Isotopies in Euclidean spaces}. (Russian) Uspehi Mat. Nauk 19 1964 no. 6 (120), 71--73.
54.78
\bibitem{dolbilin} N. P. Dolbilin,  M. A. Shtan ko, M. I.  Shtogrin.
{\it Cubic manifolds in lattices.} Izv. Ross. Akad. Nauk Ser. Mat. 58 (1994), no. 2, 93­107; translation in Russian Acad. Sci. Izv.
Math. 44 (1995), no. 2, 301­313.
\bibitem{fox1} R. H. Fox. \emph{A Quick Trip Through Knot
Theory}. Topology of 3-Manifolds and Related Topics. Prentice-Hall,
Inc., 1962.
\bibitem{FRS} Roger Fenn, Colin Rourke, Brian Sanderson. \emph{James Bundles and Applications (1996).}
http://www.maths.warwick.ac.uk/~cpr/ftp/james.ps
\bibitem{funar} L. Funar. \emph{Cubulations, immersions, mappability and a problem of Habegger.}  Ann.  scient. \'Ec. Norm. Sup.,
4e s\'erie, t. 32, 1999, pp. 681-700.
\bibitem{funar1} L. Funar. \emph{ Cubulations mod bubble moves.  Low-dimensional topology} (Funchal, 1998), Contemp. Math., 233, Amer. Math. Soc., 29-43.
\bibitem{hirsch} W. Hirsch. \emph{Smooth regular neighborhoods}. Annals
of Mathematics Vol. 76, No.3 (1962), 524-530.
\bibitem{mazur} B. Mazur. \emph{The definition of equivalence of combinatorial imbeddings.} Publications
math\'ematiques de l'I.H.\'E.S., tome 3 (1959) p.5-17.
\bibitem{matveev} S. Matveev, M. Polyak. \emph{Finite-Type Invariants of Cubic Complexes}. Acta Applicandae Mathematicae 75, pp. 125–132, 2003.
\bibitem{milnor} J. Milnor. \emph{Morse Theory}. Princeton University Press, 1969.
\bibitem{palais} R. Palais. \emph{Local triviality of the
restriction map for embeddings}. Comment. Math. Helv. 34, 1960, 305-312.
\bibitem{rolfsen} D. Rolfsen. \emph{Knots and Links}. Publish or
Perish, Inc. 1976.
\bibitem{rushing} B. Rushing. \emph{Topological
Embeddings}. Academic Press, 1973, Vol 52.
\bibitem{whitney} H. Whitney. \emph{On the topology of differentiable manifolds.}
Lectures in Topology, pp. 101-141. University of Michigan Press, Ann Arbor, Mich., 1941.
\end{thebibliography}
\end{document}